\documentclass[11pt]{article}
\usepackage{amsfonts,latexsym,rawfonts,amsmath,amssymb,amsthm,graphicx}
\numberwithin{equation}{section}
\newtheorem{theorem}{Theorem}[section]
\newtheorem{lem}[theorem]{Lemma}
\newtheorem{thm}[theorem]{Theorem}

\newtheorem{defi}[theorem]{Definition}
\newtheorem{rem}[theorem]{Remark}
\def\s{\,\,\,\,}

\def\endproof{$\hfill\Box$\\}
\def\R{\mathbb{R}}

\title{ Contact stationary Legendrian surfaces in $\mathbb{S}^5$}
\author {Yong Luo }
\date{}
\begin{document}
\maketitle
\begin{abstract}
Let $(M^5,\alpha,g_\alpha,J)$ be a 5-dimensional Sasakian  Einstein manifold with contact 1-form $\alpha$, associated metric $g_\alpha$ and almost complex structure $J$ and $L$ a contact stationary Legendrian surface in $M^5$. We will prove that $L$ satisfies the following equation
  \begin{eqnarray}\label{equ}
  -\Delta^\nu H+(K-1)H=0,
\end{eqnarray}
where $\Delta^\nu$ is the normal Laplacian w.r.t the metric $g$ on $L$ induced from $g_\alpha$ and $K$ is the Gauss curvature of $(L,g)$.

Using equation \eqref{equ} and a new Simons' type inequality for Legendrian surfaces in the standard unit sphere $\mathbb{S}^5$, we prove an integral inequality for contact stationary Legendrian surfaces in $\mathbb{S}^5$. In particular, we prove that if $L$ is a contact stationary Legendrian surface in $\mathbb{S}^5$, $B$ is the second fundamental form of $L$, $S=|B|^2$, $\rho^2=S-2H^2$ and
$$0\leq S\leq 2,$$
then we have either $\rho^2=0$ and $L$ is totally umbilic or $\rho^2\neq 0$, $S=2, H=0$ and $L$ is a flat minimal Legendrian torus.
\end{abstract}

\section{Introduction}

Let $(M^{2n+1},\alpha,g_\alpha,J)$ be a $2n+1$ dimensional contact metric manifold with contact structure $\alpha$, associated metric $g_\alpha$ and almost complex structure $J$. Assume that $(L,g)$ is an $n$-diemsional  compact Legendrian submanifold of $M^{2n+1}$ with metric $g$ induced from $g_\alpha$. The volume of $L$ is defined by
 \begin{eqnarray}
 V(L)=\int_Ld\mu,
 \end{eqnarray}
 where $d\mu$ is the volume form of $g$. A contact stationary Legendrian submanifold of $M^{2n+1}$ is a Legendrian submanifold of $M^{2n+1}$ which is a stationary point of $V$ w.r.t. Legendrian deformations. That is we call a Legendrian submanifold $L\subseteq M^{2n+1}$ a contact stationary Legendrian submanifold, if for any Legendrian deformations $L_t\subseteq M^{2n+1}$ with $L_0=L$ we have
$$\frac{dV(L_t)}{dt}|_{t=0}=0.$$
\begin{rem}
$L_t$ is a Legendrian deformation of $L:=L_0$, if $L_t$ is a Legendrian submanifold for every $t$.
\end{rem}
The E-L equation for a contact stationary Legendrian submanifold $L$ is(\cite{Ir} \cite{CLU})
\begin{eqnarray}
div_g(J H)=0,
\end{eqnarray}
where $div_g$ is the divergence w.r.t $g$ and $H$ is the mean curvature vector of $L$ in $M^{2n+1}.$
\begin{rem}
The notion of contact stationary Legendrian submanifold was first defined by  Iriyeh in \cite{Ir}  and Castro et al. in \cite{CLU} independently, where they used the name of Legendrian minimal Legendrian submanifold and contact minimal Legendrian submanifold, respectively. In this paper we prefer to use the name of contact stationary Legendrian submanifold.
\end{rem}

 The study of contact stationary Legendrian submanifolds is motivated by the study of Hamiltonian minimal Lagrangian(briefly, HSL) submanifolds, which was first studied by Ou(\cite{Oh90} \cite{Oh93}). A HSL submanifold in a K\"ahler manifold is a Lagrangian submanifold which is a stationary point of the Volume functional under Hamiltonian deformations. By \cite{Re}, Legendrian submanifolds in a Sasakian manifold $M^{2n+1}$ can be seem as links of Lagrangian submanifolds in the cone $CM^{2n+1}$, which is a K\"ahler manifold with proper metric and complex structure (see section 2). In fact, a close relation between contact stationary Legendrian submanifolds and HSL  submanifolds was found by Iriyeh \cite{Ir} and Castro et al. \cite{CLU}. Precisely, they independently proved that $C(L)$ is a HSL submanifold in $\mathbb{C}^n(n\geq2)$ if and only if $L$ is a contact stationary Legendrian submanifold in $\mathbb{S}^{2n-1}$ and $L$ is a contact stationary Legendrian submanifold in $\mathbb{S}^{2n+1}(n\geq1)$ if and only if $\Pi(L)$ is a HSL submanifold in $\mathbb{CP}^n$, where $\Pi:\mathbb{S}^{2n+1}\to \mathbb{CP}^n$ is the Hopf fibration.

From the definition we see that minimal Legendrian submaifolds are a special kind of contact stationary Legendrian submanifolds. Another special kind of contact stationary Legendrian submanifolds are Legendrian submanifolds with parallel mean curvature vector fields in the normal bundle. The study of (nonminimal)contact stationary Legendrian submanifolds of $\mathbb{S}^{2n+1}$ is relatively recent endeavor. For $n=1$, by \cite{Ir}, contact stationary Legendrian curves in $\mathbb{S}^3$ are the so called $(p,q)$ curves discovered by Schoen and Wolfson in \cite{SW}, where $p,q$ are relatively prime integers. For $n=2$, since harmonic 1-form on a 2-sphere  must be trivial, contact stationary Legendrian 2-sphere in $\mathbb{S}^5$ must be minimal and so must be the equatorial 2-spheres by Yau's result (\cite{Yau}). There are a lot of contact stationary doubly periodic surfaces form $\mathbb{R}^2$ to $\mathbb{S}^5$ by lifting H\'elein and Romon's examples (\cite{HR02}) and more contact stationary Legendrian surfaces(mainly tori) are constructed in \cite{Mi03}  \cite{Mi08}  \cite{Ir} \cite{HR05} \cite{Ma} \cite{MaS} \cite{BuC} etc.. And for general dimension examples are constructed in \cite{Oh93} \cite{Mi04} \cite{DoH} \cite{Do} \cite{Bu} \cite{JLS} \cite{Lee} \cite{CHX} etc.. See also \cite{Ono} \cite{Ka} \cite{HM} for other studies of contact stationary Legendrian submanifolds.

In this paper  we will study pinching properties of contact stationary Legendrian surfaces in $\mathbb{S}^5$. To do this we first prove an equation satisfied by contact stationary Legendrian surfaces in a Sasakian Einstein manifold, which we hope will be useful in analyzing analytic properties of contact stationary Legendrian surfaces.
\begin{thm}\label{EL equa}
Let $L$ be a contact stationary Legendrian surface in a 5-dimensional Sasakian Einstein manifold $(M^5,\alpha,g_\alpha,J)$, then $L$ satisfies the following equation:
\begin{eqnarray}\label{LeS equation}
-\Delta^\nu H+(K-1)H=0,
\end{eqnarray}
where $\Delta^\nu$ is the normal Laplacian w.r.t the metric $g$ on $L$ induced from $g_\alpha$ and $K$ is the Gauss curvature of $(L,g)$.
\end{thm}

We recall that the well-known Clifford torus is
\begin{eqnarray}\label{cliffore torus}
T_{Clif}=\mathbb{S}^1(\frac{1}{\sqrt{2}})\times \mathbb{S}^{1}(\frac{1}{\sqrt{2}})\subseteq \mathbb{S}^5.
\end{eqnarray}
In the theory of minimal surfaces, the following Simons' integral inequality and Pinching theorem due to Simons (\cite{Si}), Lawson (\cite{La}) and Chern et al. (\cite{CCK}) are well-known.
\begin{thm}[Simons, Lawson, Chern-Do Carmo-Kobayashi]\label{classic}
Let $M$ be a compact minimal surface in a unit sphere $\mathbb{S}^3$ and $B$ is the second fundamental form of $M$ in $\mathbb{S}^3$. Set $S=|B|^2$, then we have
\begin{eqnarray*}
\int_MS(2-S)d\mu\leq0.
\end{eqnarray*}
In particular, if
\begin{eqnarray*}
0\leq S\leq 2,
\end{eqnarray*}
then either $S=0$ and $M$ is totally geodesic, or $S=2$ and $M$ is the Clifford torus $T_{Clif}$, which is defined by (\ref{cliffore torus}).
\end{thm}
The above integral inequality was proved by Simons in his celebrated paper \cite{Si} and the classification result was given by Chern et al. (\cite{CCK}) and Lawson (\cite{La}), independently.

For minimal surfaces in a sphere with higher codimension, corresponding integral inequality was proved by Benko et al. (\cite{BKSS}) and Kozlowski et al. (\cite{KS}). In order to state their result, we first record an example.

\textbf{Example.} The veronese surface is a minimal surface in $\mathbb{S}^4\subseteq \mathbb{R}^5$ defined by
\begin{eqnarray*}
u: \mathbb{S}^2(\sqrt{3})\subseteq \mathbb{R}^3&\to& \mathbb{S}^4(1)\subseteq \mathbb{R}^5
\\(x,y,z)&\to& (u_1,u_2,u_3,u_4,u_5)
\end{eqnarray*}
where
\begin{eqnarray*}
u_1&=&\frac{1}{\sqrt{3}}yz, u_2=\frac{1}{\sqrt{3}}xz, u_3=\frac{1}{\sqrt{3}}xy,
\\u_4&=&\frac{1}{2\sqrt{3}}(x^2-y^2), u_5=\frac{1}{6}(x^2+y^2-2z^2).
\end{eqnarray*}
$u$ defines an isometric immersion of $\mathbb{S}^2(\sqrt{3})$ into $\mathbb{S}^4(1)$, and it maps two points $(x,y,z)$, $(-x,-y,-z)$ of $\mathbb{S}^2(\sqrt{3})$ into the same point of $\mathbb{S}^4(1)$, and so it imbeds the real projective plane into $\mathbb{S}^4(1)$.

We have
\begin{thm}[\cite{BKSS}]\label{minimal higher codimension}
Let $M$ be a minimal surface in an n-dimensional sphere $\mathbb{S}^n$, then
\begin{eqnarray}
\int_MS(2-\frac{3}{2}S)d\mu\leq0.
\end{eqnarray}
In particular, if
\begin{eqnarray*}
0\leq S\leq \frac{4}{3},
\end{eqnarray*}
then either $S=0$ and $M$ is totally geodesic, or $S=\frac{4}{3}$ ,n=4 and $M$ is the Veronese surface.
\end{thm}
The above classification for minimal surfaces in a sphere with $S=\frac{4}{3}$  was also got by Chern et al. in \cite{CCK}.

 We see that the (first) pinching constant for minimal surfaces in $\mathbb{S}^3$ is 2, but it is $\frac{4}{3}$ for minimal surfaces of higher codimensions. This is an interesting phenomenon and we think this dues to the complexity of the normal bundle, because for minimal Legendrian surfaces in $\mathbb{S}^5$, the (first) pinching constant is also 2.
 \begin{thm}[\cite{YKM}]\label{legendrian mimimal}
 If $M$ is a  minimal Legendrian surface of the unit sphere $\mathbb{S}^5$ and $0\leq S\leq 2$, then $S$ is identically 0 or 2.
 \end{thm}

\begin{rem}
For higher dimensional case of this theorem we refer to \cite{DV}.
\end{rem}

 All of these results are based on calculating the Laplacian of $S$ and then get Simons' type equalities or inequalities,  a powerful method which was originated from \cite{Si}. The minimal condition is used to cancel some terms in the resulting calculation and to some extent it is important. In this note we prove a Simons' type inequality (lemma \ref{main result}) for Legendrian surfaces in $\mathbb{S}^5$, without minimal condition. By using equation (\ref{LeS equation}) and this Simons' type inequality we get
\begin{thm}\label{inte ineq}
Let $L:\Sigma\to \mathbb{S}^5$ be a contact stationary Legendrian surface, where $\mathbb{S}^5$ is the unit sphere with standard contact structure and metric (as given in the end of section 2). Then we have
\begin{eqnarray*}
\int_L\rho^2(3-\frac{3}{2}S+2H^2)d\mu\leq0,
\end{eqnarray*}
where $\rho^2:=S-2H^2$. In particular, if
\begin{eqnarray*}
0\leq S\leq 2,
\end{eqnarray*}
then either $\rho^2=0$ and $L$ is totally umbilic, or $\rho^2\neq 0$, $S=2, H=0$ and $L$ is a flat minimal Legendrian torus.
\end{thm}
\begin{rem}
Because minimal Legendrian surfaces are contact stationary Legendrian surfaces and for minimal Legendrian surfaces $\rho^2=S$ and totally umbilic minimal surfaces are totally geodesic, we see that theorem \ref{legendrian mimimal} is a corollary of theorem \ref{inte ineq}.
\end{rem}

Integral inequality and gap phenomenon for submanifolds satisfying a fourth order quasi-elliptic nonlinear equation was first studied by Li. In \cite{Li1} \cite{Li2} and \cite{Li02}, Li proved several gap theorems for Willmore submanifolds in a sphere. These results are partial motivations of our paper.

We end this introduction by recalling a classification theorem of flat minimal Legendrian toruses in $\mathbb{S}^5$. For a constant $\theta$ let $T_\theta$ be the 2-torus in $\mathbb{S}^5$ defined by
$$T_\theta=\{(z_1,z_2,z_3)\in \mathbb{C}^3:|z_i|=\frac{1}{3}, i=1,2,3 \s and \s \sum_iargz_i=\theta\}.$$
$T_\theta$ is called the generalized Clifford torus and it is a flat minimal Legendrian torus in $\mathbb{S}^5$. Its projection under the Hopf map $\pi:\mathbb{S}^5\to \mathbb{CP}^2$ is a flat minimal Lagrangian torus , which is also called a generalized Clifford torus. It is proved in \cite{LOY} that a flat minimal Lagrangian torus in $\mathbb{CP}^2$ must be $\mathbb{S}^1\times \mathbb{S}^1$. By the correspondence of minimal Lagrangian surfaces in $\mathbb{CP}^2$ and minimal Legendrian surfaces in $\mathbb{S}^5$ (cf.\cite{Re}), we see that a flat minimal Legendrian torus in $\mathbb{S}^5$ must be a generalized Clifford torus. For more details we refer to \cite{Ha}, page 853.

The rest of this paper is organized as follows: In section 2 we collect some basic material from Sasakian geometry, which will be used in the next section. In section 3 we prove our main results, theorem\ref{EL equa}, and theorem \ref{inte ineq}.

\section{Preliminaries on contact geometry}
In this section we recall some basic material from contact geometry. For more information we refer to \cite{Bl}.
\subsection{Contact Manifolds}
\begin{defi}
A contact manifold $M$ is an odd dimensional manifold with a one form $\alpha$ such that $\alpha\wedge(d\alpha)^n\neq0$, where $dimM=2n+1$.
\end{defi}
Assume now that $(M,\alpha)$ is a given contact manifold of dimension $2n+1$. Then $\alpha$ defines a $2n$-dimensional vector bundle over $M$, where the fibre at each point $p\in M$ is given by
$$\xi_p=Ker\alpha_p.$$
Sine $\alpha\wedge (d\alpha)^n$ defines a volume form on $M$, we see that
$$\omega:=d\alpha$$
is a closed nondegenerate 2-form on $\xi\oplus\xi$ and hence it defines a symplectic product on $\xi$ such that $(\xi,\omega|_{\xi\oplus\xi})$ becomes a symplectic vector bundle. A consequence of this fact is that there exists an almost complex bundle structure
$$\tilde{J}:\xi\to\xi$$
compatible with $d\alpha$, i.e. a bundle endomorphism satisfying:
\\(1) $\tilde{J}^2=-id_\xi$,
\\(2) $d\alpha(\tilde{J}X,\tilde{J}Y)=d\alpha(X,Y)$ for all $X,Y\in\xi$,
\\(3) $d\alpha(X,\tilde{J}X)>0$ for $X\in\xi\setminus {0}$.

Since $M$ is an odd dimensional manifold, $\omega$ must be degenerate on $TM$, and so we obtains a line bundle $\eta$ over $M$ with fibres
$$\eta_p:=\{V\in T_pM|\omega(V,W)=0 \s\forall\s W\in\xi_p\}.$$
\begin{defi}
The Reeb vector field $\textbf{R}$ is the section of $\eta$ such that $\alpha(\textbf{R})=1$.
\end{defi}

Thus $\alpha$ defines a splitting of $TM$ into a line bundle $\eta$ with the canonical section $\textbf{R}$ and a symplectic vector bundle $(\xi,\omega|\xi\oplus\xi)$. We denote the projection along $\eta$ by $\pi$, i.e.
\begin{eqnarray*}
&&\pi:TM\to\xi,
\\&&\pi(V):=V-\alpha(V)\textbf{R}.
\end{eqnarray*}
Using this projection we extend the almost complex structure $\tilde{J}$ to a section $J\in\Gamma(T^*M\otimes TM)$ by setting
$$J(V)=\tilde{J}(\pi(V)),$$
for $V\in TM$.

We call $J$ an almost complex structure of the contact manifold $M$.
\begin{defi}
Let $(M,\alpha)$ be a contact manifold, a submanifold $L$ of $(M,\alpha)$ is called an isotropic submanifold if $T_xL\subseteq\xi_x$ for all $x\in L$.
\end{defi}
For algebraic reasons the dimension of an isotropic submanifold of a $2n+1$ dimensional contact manifold can not be bigger than $n$.
\begin{defi}
An isotropic submanifold $L\subseteq(M,\alpha)$ of maximal possible dimension $n$ is called a Legendrian submanifold.
\end{defi}
\subsection{Sasakian manifolds}
Let $(M,\alpha)$ be a contact manifold, with the almost complex structure $J$ and Reeb field $\textbf{R}$. A Riemannian metric $g_\alpha$ defined on $M$ is said to be associated, if it satisfies the following three conditions:
\\(1) $g_\alpha(\textbf{R},\textbf{R})=1$,
\\(2) $g_\alpha(V,\textbf{R})=0$, $\forall\s V\in\xi$,
\\(3) $\omega(V,JW)=g_\alpha(V,W)$, $\forall\s V,W\in\xi$.\

We should mention here that on any contact manifold there exists an associated metric on it, because we can construct one in the following way. We introduce a bilinear form $b$ by
$$b(V,W):=\omega(V,JW),$$
then the tensor
$$g:=b+\alpha\otimes\alpha$$
defines an associated metric on $M$.

Sasakian manifolds are the odd dimensional analogue of K\"ahler manifolds. They are defined as follows.
\begin{defi}
A contact manifold $(M,\alpha)$ with an associated metric $g_\alpha$ is called Sasakian, if the cone $CM$ equipped with the following extended metric $\bar{g}$
\begin{eqnarray}\label{cone metic}
(CM,\bar{g})=(\mathbb{R}_+\times M,dr^2+r^2g_\alpha)
\end{eqnarray}
is K\"ahler w.r.t the following canonical almost complex structure $J$ on $TCM=\mathbb{R}\oplus\langle\textbf{R}\rangle\oplus\xi:$
$$J(r\partial r)=\textbf{R}, J(\textbf{R})=-r\partial r.$$
Furthermore if $g_\alpha$ is Einstein, $M$ is called a Sasakian  Einstein manifold.
\end{defi}
We record several lemmas which are well known in Sasakian geometry. These lemmas will be used in the next section.
\begin{lem}
Let $(M,\alpha,g_\alpha,J)$ be a Sasakian manifold. Then
\begin{eqnarray}\label{Reeb}
\bar{\nabla}_X\textbf{R}=-JX,
\end{eqnarray}
and
\begin{eqnarray}\label{derivatives}
(\bar{\nabla}_XJ)(Y)=g(X,Y)\textbf{R}-\alpha(Y)X,
\end{eqnarray}
for $X,Y\in TM$, where $\bar{\nabla}$ is the Levi-Civita connection on $(M,g_\alpha)$.
\end{lem}
\begin{lem}\label{mean curvatue form}
Let $L$ be a Legendrian submanifold in a Sasakian Einstein manifold $(M,\alpha,g_\alpha,J)$, then the mean curvature form $\omega(H,\cdot)|_L$ defines a closed one form on $L$.
\end{lem}
For a proof of this lemma we refer to \cite{Le}, Proposition A.2 or \cite{Sm}, lemma 2.8. In fact they proved this result under a weaker assumption that $(M,\alpha,g_\alpha,J)$ is a weakly Sasakian Einstein manifold, where weakly Einstein means that $g_\alpha$ is Einstein only when restricted to the contact hyperplane $Ker\alpha$.
\begin{lem}\label{orthogonal}
Let $L$ be a Legendrian submanifold in a Sasakian manifold $(M,\alpha,g_\alpha,J)$ and $B$ be the second fundamental form of $L$ in $M$. Then we have
\begin{eqnarray}
g_\alpha(B(X,Y),\textbf{R})=0,
\end{eqnarray}
for any $X,Y\in TL$.
\end{lem}
\proof For any $X,Y\in TL$,
\begin{eqnarray*}
\langle B(X,Y),\textbf{R}\rangle&=&\langle\bar{\nabla}_XY,\textbf{R}\rangle
\\&=&-\langle Y,\bar{\nabla}_X\textbf{R}\rangle
\\&=&\langle Y,JX\rangle
\\&=&\omega(X,Y)
\\&=&d\alpha(X,Y)
\\&=&0,
\end{eqnarray*}
where in the third equality we used (\ref{Reeb}).
\endproof

In particular this lemma implies that the mean curvature $H$ of $L$ is orthogonal to the Reeb field $\textbf{R}$.
\begin{lem}\label{commute of J}
For any $Y,Z\in Ker\alpha$, we have
\begin{eqnarray}
g_\alpha(\bar{\nabla}_X(JY),Z)=g_\alpha(J\bar{\nabla}_XY,Z).
\end{eqnarray}
\end{lem}
\proof Note that
$$(\bar{\nabla}_XJ)Y=\bar{\nabla}_X(JY)-J\bar{\nabla}_XY.$$
Therefore by using (\ref{derivatives}) we have
\begin{eqnarray*}
\langle \bar{\nabla}_X(JY),Z\rangle&=&\langle(\bar{\nabla}_XJ)Y,Z \rangle+\langle J\bar{\nabla}_XY,Z\rangle
\\&=&\langle J\bar{\nabla}_XY,Z\rangle,
\end{eqnarray*}
for any $Y,Z\in Ker\alpha$. \endproof

A most canonical example of Sasakian Einstein manifolds is the standard odd dimensional sphere $\mathbb{S}^{2n+1}$.

\textbf{The standard sphere $\mathbb{S}^{2n+1}$.}
 Let $\mathbb{C}^n=\mathbb{R}^{2n+2}$ be the Euclidean space with coordinates $(x_1,...,x_{n+1},y_1,...,y_{n+1})$ and $\mathbb{S}^{2n+1}$ be the standard unit sphere in $\mathbb{R}^{2n+2}$. Define
$$\alpha_0=\frac{1}{2}\sum_{j+1}^{n+1}(x_jdy_j-y_jdx_j),$$
then
$$\alpha:=\alpha_0|_{\mathbb{S}^{2n+1}}$$
defines a contact one form on $\mathbb{S}^{2n+1}$. Assume that $g_0$ is the standard metric on $\mathbb{R}^{2n+2}$ and $J_0$ is the standard complex structure of $\mathbb{C}^n$. We define
$$g_\alpha=g_0|_{\mathbb{S}^{2n+1}}, J=J_0|_{\mathbb{S}^{2n+1}},$$
then $(\mathbb{S}^{2n+1},\alpha,g_\alpha,J)$ is a Sasakian Einstein manifold with associated metric $g_\alpha$. Its contact hyperplane is characterized by
$$Ker\alpha_x=\{Y\in T_x\mathbb{S}^{2n+1}|\langle Y,Jx\rangle=0\}.$$

\section{Proof of the theorems}
\subsection{Several lemmas}
In this part we assume that $(M,\alpha,g_\alpha,J)$ is a Sasakian manifold. We show several lemmas which are analogous results in K\"ahler geometry.

The first lemma shows $\omega=d\alpha$ when restricted to the contact hyperplane $Ker\alpha$ behaviors as the K\"ahler form on a K\"ahler manifold.
\begin{lem}\label{kahler for w}
Let $X,Y,Z\in Ker\alpha$, then
\begin{eqnarray}
\bar{\nabla}_X\omega(Y,Z)=0,
\end{eqnarray}
where $\bar{\nabla}$ is the derivative w.r.t $g_\alpha$.
\end{lem}
\proof
\begin{eqnarray*}
\bar{\nabla}_X\omega(Y,Z)&=&X(\omega(Y,Z))-\omega(\bar{\nabla}_XY,Z)-\omega(Y,\bar{\nabla}_XZ)
\\&=&-Xg_\alpha(Y,JZ)-\omega(\bar{\nabla}_XY,Z)-\omega(Y,\bar{\nabla}_XZ)
\\&=&-g_\alpha(\bar{\nabla}_XY,JZ)-g_\alpha(Y,\bar{\nabla}_XJZ)
+g_\alpha(\bar{\nabla}_XY,JZ)+g_\alpha(Y,J\bar{\nabla}_XZ)
\\&=&0,
\end{eqnarray*}
where in the third equality we used $g_\alpha(Y,\bar{\nabla}_XJZ)=g_\alpha(Y,J\bar{\nabla}_XZ)$, which is a direct corollary of (\ref{derivatives}).
\endproof

Now let $L$ be a Legendrian submanifold of $M$. We have a natural identification of $NL\cap Ker\alpha$ with $T^\ast L$, where $NL$ is the normal bundle of $L$ and $T^\ast L$ is the cotangent bundle.
\begin{defi}
$\tilde{\omega}:NL\cap Ker\alpha\to T^\ast L$ is the bundle isomorphism defined by
$$\tilde{\omega}_p(v_p)=(v_p\rfloor \omega_p)|_{T_pL},$$
where $p\in L$ and $v_p\in (NL\cap Ker\alpha)_p.$
\end{defi}
Recall that $\omega(\textbf{R})=0$ and $g_\alpha(V,W)=\omega(V,JW)$ for any $V,W\in\xi$, hence $\tilde{\omega}$ defines an isomorphism.

We have
\begin{lem}\label{main lem1}
Let $V\in\Gamma (NL\cap Ker\alpha)$. Then
\begin{eqnarray}\label{main equ}
\tilde{\omega}(\Delta^\nu V-\langle\Delta^\nu V,\textbf{R}\rangle\textbf{R}+V)&=&\Delta(\tilde{\omega}(V))\s i.e. \nonumber
\\(\Delta^\nu V+V)\rfloor\omega&=&\Delta(V\rfloor\omega),
\end{eqnarray}
where $\Delta$ is the Laplace-Beltrami operator on $(L,g)$.
\end{lem}
\begin{rem}
This kind of lemma in the context of symplectic geometry was proved by Oh (\cite{Oh90}, lemma 3.3). Our proof follows his argument with only slight modifications.
\end{rem}
\proof We first show that
\begin{eqnarray}\label{important equality}
\nabla_X(\tilde{\omega}(V))=\tilde{\omega}(\nabla_X^\nu V-\langle\nabla_X^\nu V,\textbf{R}\rangle\textbf{R})
\end{eqnarray}
for any $X\in TL$.
Equality (\ref{important equality}) is equivalent to
\begin{eqnarray}\label{1}
\nabla_X(\tilde{\omega}(V))(Y)=\tilde{\omega}(\nabla_X^\nu V-\langle\nabla_X^\nu V,\textbf{R}\rangle\textbf{R})(Y)
\end{eqnarray}
for any $Y\in TL$.
\begin{eqnarray*}
\nabla_X(\tilde{\omega}(V))(Y)&=&\nabla_X(\tilde{\omega}(V)(Y))-\tilde{\omega}(V)(\nabla_XY)
\\&=&\bar{\nabla}_X(\omega(V,Y))-\tilde{\omega}(V)(\nabla_XY)
\\&=&\omega(\nabla^\nu_XV,Y)+\omega(V,\nabla_XY)-\omega(V,\nabla_XY)
\\&=&\omega(\nabla^\nu_XV,Y)
\\&=&\tilde{\omega}(\nabla^\nu_XV-\langle\nabla_X^\nu V,\textbf{R}\rangle\textbf{R})(Y),
\end{eqnarray*}
where in the third equality we used $\bar{\nabla}_X\omega=0,$ when restricted to $Ker\alpha$, which is proved in lemma \ref{kahler for w}.

Let $p\in L$ and we choose an orthonormal frame $\{E_1,...,E_n\}$ on $TL$ with $\nabla_{E_i}E_j(p)=0$, then the general Laplacian $\Delta$ can be written as
$$\Delta\psi(p)=\sum_{i=1}^n\nabla_{E_i}\nabla_{E_i}\psi(p),$$
where $\psi$ is a tensor on $L$. Therefore
\begin{eqnarray*}
&&(\tilde{\omega}^{-1}\circ\Delta\cdot\tilde{\omega}(V))(p)
\\&=&(\tilde{\omega}^{-1}\circ\sum_{i=1}^n\nabla_{E_i}\nabla_{E_i}\tilde{\omega}(V))(p)
\\&=&\sum_{i=1}^n(\tilde{\omega}^{-1}\nabla_{E_i}\tilde{\omega}
\cdot\tilde{\omega}^{-1}\nabla_{E_i}\tilde{\omega}(V))(p)
\\&=&\sum_{i=1}^n(\tilde{\omega}^{-1}\nabla_{E_i}\tilde{\omega}(\nabla_{E_i}^\nu V-\langle\nabla_{E_i}^\nu V,\textbf{R}\rangle\textbf{R})(p)
\\&=&\sum_{i=1}^n\nabla_{E_i}^\nu(\nabla_{E_i}^\nu V-\langle\nabla_{E_i}^\nu V,\textbf{R}\rangle\textbf{R})-\langle\nabla_{E_i}^\nu(\nabla_{E_i}^\nu V-\langle\nabla_{E_i}^\nu V,\textbf{R}\rangle\textbf{R}),\textbf{R}\rangle\textbf{R}
\\&=&\Delta^\nu V-\langle\Delta^\nu V,\textbf{R}\rangle\textbf{R}+V,
\end{eqnarray*}
where in the third and fourth equalities we used (\ref{important equality}) and in the last equality we used equality (\ref{Reeb}).
\endproof
\subsection{Proof of theorem \ref{EL equa}}
We see that for any function $s$ defined on $L$,
\begin{eqnarray*}
0&=&\int_LsdivJH d\mu
=\int_Lg(J H,\nabla s)d\mu
\\&=&\int_L\omega(H,\nabla s)d\mu
=\int_L\langle \omega\rfloor H,\omega \rfloor\nabla s\rangle d\mu
\\&=&\int_L\langle\omega\rfloor H,d s \rangle
=\int_L\delta(\omega\rfloor H)sd\mu.
\end{eqnarray*}
Therefore the E-L equation for $L$ is equivalent to
\begin{eqnarray}
\delta(\omega\rfloor H)=0,
\end{eqnarray}
where $\delta$ is the adjoint operator of $d$ on $L$.

By lemma \ref{mean curvatue form} we see that $L$ satisfies
\begin{eqnarray}
\Delta_h(\omega\rfloor H)=0,
\end{eqnarray}
where $\Delta_h:=\delta d+d\delta$ is the Hodge-Laplace operator. That is the mean curvature form of $L$ is a harmonic one form.

To proceed on, we need the following Weitzenb\"ock formula
\begin{lem}
Let $M$ be an $n$ dimensional oriented Riemannian manifold. If $\{V_i\}$ is a local orthonormal frame field and $\{\omega^i\}$ is its dual co-frame field, then
$$\Delta_h=-\sum_iD^2_{V_iV_i}+\sum_{ij}\omega^i\wedge i(V_j)R_{V_iV_j},$$
where $D^2_{XY}\equiv D_XD_Y-D_{D_XY}$ represents the covariant derivatives, $\Delta_d=d\delta+\delta d$ is the Hodge-Laplace and $R_{XY}=-D_XD_Y+D_YD_X+D_{[X,Y]}$ is the curvature tensor.
\end{lem}
\begin{rem}
For a detailed discussion on the Weitzenb\"ock formula we refer to Wu (\cite{Wu}).
\end{rem}

Using the Weitzenb\"ock formula we have
\begin{eqnarray}\label{import 2}
-\Delta(\omega\rfloor H)+\sum_{ij}\omega^i\wedge i(V_j)R_{V_iV_j}(\omega\rfloor H)=0,
\end{eqnarray}
where $\{V_i\}$ is a local orthogonal frame field and $\{\omega^i\}$ is its dual co-frame field on $L$.

 Denote $\omega\rfloor H$ by $\theta_H=\sum_k\theta_k\omega^k$, we have
\begin{eqnarray*}
\sum_{ij}\omega^i\wedge i(V_j)R_{V_iV_j}\theta_H
 &=&\sum_{ij}R_{V_iV_j}\theta_H(V_j)\omega^i
 \\&=&\sum_{ijk}R_{V_iV_j}\omega^k(V_j)\theta_k\omega^i
 \\&=&-\sum_{ijk}\omega^k(R_{V_iV_j}V_j)\theta_k\omega^i
 \\&=&-\sum_{ijk}\langle R_{V_iV_j}V_j, V_k\rangle\theta_k\omega^i
\\ &=&-\sum_{ij}\langle R_{V_iV_j}V_j, V_i\rangle\theta_i\omega^i
 \\&=&K\theta_H.
\end{eqnarray*}
That is
\begin{eqnarray}\label{important 3}
\sum_{ij}\omega^i\wedge i(V_j)R_{V_iV_j}(\omega\rfloor H)=K\omega\rfloor H.
\end{eqnarray}

Recall that $H\in NL\cap Ker\alpha$, using (\ref{main equ}) to $H$ we get
\begin{eqnarray}\label{important equality 2}
\Delta(\omega\rfloor H)=(\Delta^\nu H+H)\rfloor\omega.
\end{eqnarray}
Combining (\ref{import 2})-(\ref{important equality 2}), we have
\begin{eqnarray*}
0&=&-\Delta^\nu H\rfloor\omega-H+K\omega\rfloor H
\\&=&(-\Delta^\nu H+(K-1)H)\rfloor\omega,
\end{eqnarray*}
which implies that
\begin{eqnarray}
-\Delta^\nu H+(K-1)H=f\textbf{R}
\end{eqnarray}
for some function $f$ on $L$.

The next lemma is one of our key observations which states that a Legendrian submanifold in a Sasakian manifold is contact stationary if and only if $\langle\Delta^\nu H,\textbf{R}\rangle=0$.
\begin{lem}
Let $L\subseteq (M^{2n+1},\alpha,g_\alpha,J)$ be a contact stationary Legendrian submanifold. Then we have $\Delta^\nu H$ is orthogonal to $\textbf{R}$.
\end{lem}
\proof For any point $p\in L$, we choose a local orthonormal frame $\{E_i:i=1,...,n\}$ of $L$ such that $\nabla_{E_i}E_j(p)=0$. We have at $p$ (in the following computation we adopt the Einstein summation convention)
\begin{eqnarray*}
\langle\Delta^\nu H,\textbf{R}\rangle&=&\langle\nabla^\nu_{E_i}\nabla^\nu_{E_i}H,\textbf{R}\rangle
\\&=&E_i\langle\nabla^\nu_{E_i}H,\textbf{R}\rangle
-\langle\nabla^\nu_{E_i}H,\bar{\nabla}_{E_i}\textbf{R}\rangle
\\&=&E_i\langle\nabla^\nu_{E_i}H,\textbf{R}\rangle+\langle\nabla^\nu_{E_i}H,JE_i\rangle
\\&=&E_i(E_i\langle H,\textbf{R}\rangle-\langle H,\bar{\nabla}_{E_i}\textbf{R}\rangle)+\langle\nabla^\nu_{E_i}H,JE_i\rangle
\\&=&E_i\langle H,JE_i\rangle+\langle\nabla^\nu_{E_i}H,JE_i\rangle
\\&=&2\langle\nabla^\nu_{E_i}H,JE_i\rangle+\langle H,\bar{\nabla}_{E_i}JE_i\rangle
\\&=&2\langle\nabla^\nu_{E_i}H,JE_i\rangle+\langle H,J\bar{\nabla}_{E_i}E_i\rangle
\\&=&2\langle\nabla^\nu_{E_i}H,JE_i\rangle
\\&=&2\langle\bar{\nabla}_{E_i}H,JE_i\rangle
\\&=&-2\langle J\bar{\nabla}_{E_i}H,E_i\rangle
\\&=&-2\langle\bar{\nabla}_{E_i}JH,E_i\rangle
\\&=&-2\langle\nabla_{E_i}JH,E_i\rangle
\\&=&-2div_g(JH)
\\&=&0.
\end{eqnarray*}
Note that in this computation we used lemma \ref{derivatives}, lemma \ref{orthogonal} and lemma \ref{commute of J} several times and the last equality holds because $L$ is contact stationary. \endproof

Therefore we have $$(-\Delta^\nu H+(K-1)H)\bot\textbf{R}$$ by this lemma and lemma \ref{orthogonal}, which shows $f\equiv0$, i.e.
$$-\Delta^\nu H+(K-1)H=0,$$
and we are done.
\endproof
\subsection{Proof of theorem \ref{inte ineq}}
Let $L$ be a Legendrian surface in $\mathbb{S}^5$ with the induced metric $g$. Let $\{e_1,e_2\}$ be an orthogonal frame on $L$ such that $\{e_1,e_2,Je_1,Je_2,\textbf{R}\}$ be a orthonormal frame on $\mathbb{S}^5$.

In the following we use indexes $i,j,k,l,s,t,m$ and $\beta,\gamma$ such that
\begin{eqnarray*}
1\leq i,j,k,l,s,t,m&\leq&2,
\\1\leq\beta,\gamma&\leq&3,
\\ \gamma^\ast=\gamma+2,\s \beta^\ast&=&\beta+2.
\end{eqnarray*}

Let $B$ be the second fundamental form of $L$ in  $\mathbb{S}^5$ and define
\begin{eqnarray}
h_{ij}^k&=&g_\alpha(B(e_i,e_j),Je_k),
\\h^3_{ij}&=&g_\alpha(B(e_i,e_j),\textbf{R}).
\end{eqnarray}
Then
\begin{eqnarray}
h_{ij}^k&=&h_{ik}^j=h_{kj}^i,
\\h^3_{ij}&=&0.
\end{eqnarray}
The Gauss equations and Ricci equations are
\begin{eqnarray}
R_{ijkl}&=&(\delta_{ik}\delta_{jl}-\delta_{il}\delta_{jk})+\sum_s(h^s_{ik}h^s_{jl}-h^s_{il}h^s_{jk})\label{basic equation 1}
\\R_{ik}&=&\delta_{ik}+2\sum_sH^sh^s_{ik}-\sum_{s,j}h^s_{ij}h^s_{jk},
\\2K&=&2+4H^2-S,
\\R_{3412}&=&\sum_i(h_{i1}^1h_{i2}^2-h_{i2}^1h_{i1}^2)\nonumber
\\&=&\det h^1+\det h^2,
\end{eqnarray}
where $h^1,h^2$ are the second fundamental forms w.r.t. the directions $Je_1$, $Je_2$ respectively.

 In addition we have the following Codazzi equations and Ricci identities
\begin{eqnarray}
h^\beta_{ijk}&=&h^\beta_{ikj},
\\h^\beta_{ijkl}-h^\beta_{ijlk}&=&\sum_mh^\beta_{mj}R_{mikl}+\sum_mh^\beta_{mi}R_{mjkl}+\sum_\gamma h^\gamma_{ij}R_{\gamma^\ast\beta^\ast kl}.\label{basic equation 2}
\end{eqnarray}

Using these equations, we can get the following Simons' type inequality:
\begin{lem}\label{main result}
Let $L$ be a Legendrian surface in $\mathbb{S}^5$. Then we have
\begin{eqnarray}\label{main lemma}
\frac{1}{2}\Delta\sum_{i,j,\beta}(h^\beta_{ij})^2&\geq&|\nabla^T h|^2-2|\nabla^T H|^2-2|\nabla^\nu H|^2 +\sum_{i,j,k,\beta}(h^\beta_{ij}h^\beta_{kki})_j \nonumber
\\&+&S-2H^2+2(1+H^2)\rho^2-\rho^4-\frac{1}{2}S^2,
\end{eqnarray}
where $|\nabla^T h|^2=\sum_{i,j,k,s}(h^s_{ijk})^2$ and $|\nabla^T H|^2=\sum_{i,s}(H^s_i)^2$.
\end{lem}
\proof Using equations from (\ref{basic equation 1}) to (\ref{basic equation 2}), we have
\begin{eqnarray}\label{simon type}
\frac{1}{2}\Delta\sum_{i,j,\beta}(h^\beta_{ij})^2
&=&\sum_{i,j,k,\beta}(h^\beta_{ijk})^2+\sum_{i,j,k,\beta}h^\beta_{ij}h^\beta_{kijk}\nonumber
\\&=&|\nabla h|^2-4|\nabla^\nu H|^2+\sum_{i,j,k,\beta}(h^\beta_{ij}h^\beta_{kki})_j+\sum_{i,j,l,k,\beta} h^\beta_{ij}(h^\beta_{lk}R_{lijk}+h^\beta_{il}R_{lj})\nonumber
\\&+&\sum_{i,j,k,\beta,\gamma} h^\beta_{ij}h^\gamma_{ki}R_{\gamma^\ast\beta^\ast jk}\nonumber
\\&=&|\nabla h|^2-4|\nabla^\nu H|^2+\sum_{i,j,k,s}(h^s_{ij}h^s_{kki})_j+2K\rho^2-2(\det h^1+\det h^2)^2\nonumber
\\&\geq&|\nabla h|^2-4|\nabla^\nu H|^2+\sum_{i,j,k,\beta}(h^\beta_{ij}h^\beta_{kki})_j+2(1+H^2)\rho^2-\rho^4-\frac{1}{2}S^2,
\end{eqnarray}
where $\rho^2:=S-2H^2$ and in the above calculations we used the following identities
\begin{eqnarray*}
\sum_{i,j,k,l,\beta} h^\beta_{ij}(h^\beta_{lk}R_{lijk}+h^\beta_{il}R_{lj})&=&2K\rho^2,
\\\sum_{i,j,k,\beta,\gamma} h^\beta_{ij}h^\gamma_{ki}R_{\gamma^\ast\beta^\ast jk}&=&-2(\det h^1+\det h^2)^2,
\end{eqnarray*}
where in the first equality we used $R_{lijk}=K(\delta_{lj}\delta_{ik}-\delta_{lk}\delta_{ij})$ and $R_{lj}=K\delta_{lj}$ in a proper coordinate, because $L$ is a surface.

Note that
\begin{eqnarray}\label{main idea1}
|\nabla h|^2&=&\sum_{i,j,k,\beta}(h^\beta_{ijk})^2\nonumber
\\&=&|\nabla^T h|^2+\sum_{i,j,k}(h^3_{ijk})^2\nonumber
\\&=&|\nabla^T h|^2+\sum_{i,j,k}(h^k_{ij})^2\nonumber
\\&=&|\nabla^T h|^2+S,
\end{eqnarray}
where in the third equality we used
\begin{eqnarray*}
h^3_{ijk}&=&\langle\bar{\nabla}_{e_k}B(e_i,e_j),\textbf{R}\rangle
\\&=&-\langle B(e_i,e_j),\bar{\nabla}_{e_k}\textbf{R}\rangle
\\&=&\langle B(e_i,e_j),Je_k\rangle
\\&=&h^k_{ij}.
\end{eqnarray*}
Similarly we have
\begin{eqnarray}\label{main idea2}
|\nabla^\nu H|^2=|\nabla^TH|^2+H^2.
\end{eqnarray}
Combing (\ref{simon type}), (\ref{main idea1}) and (\ref{main idea2}) we get (\ref{main lemma}).
\endproof

Now we prove an integral equality for $L$, by using the equation (\ref{LeS equation}).
\begin{lem}
Let $L:\Sigma\to \mathbb{S}^5$ be a contact stationary Legendrian surface, where $\mathbb{S}^5$ is the unit sphere with standard contact structure and metric. Then
\begin{eqnarray}\label{integral equality}
\int_L|\nabla^\nu H|^2d\mu=-\int_L(K-1)H^2d\mu,
\end{eqnarray}
where $|\nabla^\nu H|^2=\sum_{\beta,i}(H^\beta_i)^2$.
\end{lem}
\proof By using (\ref{LeS equation}) we have
\begin{eqnarray}\label{important 4}
|\nabla^\nu H|^2&=&\sum_{\beta,i}(H^\beta_i)^2 \nonumber
\\&=&\sum_{\beta,i}(H^\beta_iH^\beta)_i-\sum_\beta H^\beta\Delta^\nu H^\beta \nonumber
\\&=&\sum_{\beta,i}(H^\beta_iH^\beta)_i-(K-1)H^2.
\end{eqnarray}
We get (\ref{integral equality}) by integrating over (\ref{important 4}).
\endproof

Integrating over (\ref{main lemma}) and using $|\nabla^Th|^2\geq 3|\nabla^TH|^2$ (see appendix, Lemma \ref{appendix}) we get
\begin{eqnarray}\label{5.11}
0&\geq&\int_L[(|\nabla^T h|^2-2|\nabla^T H|^2)-2|\nabla^\nu H|^2+S-2H^2+2(1+H^2)\rho^2-\rho^4-\frac{1}{2}S^2]d\mu \nonumber
\\ &\geq& \int_L[-2|\nabla^\nu H|^2+S-2H^2+2(1+H^2)\rho^2-\rho^4-\frac{1}{2}S^2]d\mu \nonumber
\\&=&\int_L(2-\rho^2)\rho^2d\mu+\int_L 2H^2\rho^2+2(K-1)H^2-2H^2+S-\frac{1}{2}S^2d\mu \nonumber
\\&=&\int_L(2-\rho^2)\rho^2d\mu+\int_L 2H^2\rho^2+(4H^2-S)H^2-2H^2+S-\frac{1}{2}S^2d\mu \nonumber
\\&=&\int_L(2-\rho^2)\rho^2d\mu+\int_LH^2S-2H^2+S-\frac{1}{2}S^2d\mu \nonumber
\\&=&\int_L(2-\rho^2)\rho^2d\mu+\int_LH^2(S-2)+\frac{S}{2}(2-S)d\mu\nonumber
\\&=&\int_L(2-\rho^2)\rho^2+(2-S)(\frac{S}{2}-H^2)d\mu\nonumber
\\ &=& \int_L\rho^2(2-\rho^2)+\frac{\rho^2}{2}(2-S)d\mu \nonumber
\\&=&\int_L\frac{3}{2}\rho^2(2-S)+2H^2\rho^2d\mu,\nonumber
\end{eqnarray}
where in the second equality we used the Gauss equation $2K=2+4H^2-S$.

Therefore we obtain the desired integral inequality
$$\int_L\rho^2(3-\frac{3}{2}S+2H^2)d\mu\leq0.$$
Particularly if $0\leq S\leq2$, we must have $\rho^2=0$ and $L$ is totally umbilic or  $\rho^2\neq 0$, which implies $S=2, H=0$ and $L$ is a flat minimal Legendrian torus.
\endproof

\section{Appendix}
In this section we prove the following lemma.
\begin{lem}\label{appendix}
Let $L$ be a Legendrian surface in $\mathbb{S}^5$, and assume that  $|\nabla^T h|^2$, $|\nabla^T H|^2$ are defined in Lemma \ref{main result}. Then we have
$$|\nabla^T h|^2\geq 3|\nabla^T H|^2.$$
\end{lem}
\proof  We construct the flowing symmetric tracefree tensor:
\begin{eqnarray}
F_{ijk}^s=h_{ijk}^s-\frac{1}{2}(H_{i}^s\delta_{jk}+H_{j}^s\delta_{ik}+H_{k}^s\delta_{ji}).
\end{eqnarray}
Then it is easy to see that
$$|F|^2=|\nabla^T h|^2-3|\nabla^T H|^2$$
and we get $|\nabla^T h|^2\geq 3|\nabla^T H|^2$. \endproof

\textbf{Final discussions.} At the end of this paper we propose several questions which we will study in the future.
\\\textbf{ Problem 1:} Is any umbilical contact stationary Legendrian surface in $\mathbb{S}^5$ with $0\leq S\leq 2$ totally geodesic?
\\\textbf{Problem 2:} Assume that $L$ is a closed csL submanifold in $\mathbb{S}^{2n+1}$, satisfying $0\leq S\leq n$, then is $L$ totally geodesic or $S=n$?
\\\textbf{Problem 3:} Is any contact stationary Legendrian surface in $\mathbb{S}^5$ with second fundamental form of constant length minimal?
\\\textbf{Problem 4:} What is the second gap for minimal Legendrian submanifolds in a sphere?

\vspace{1cm}

\textbf{Acknowledgement.} I would like to thank  professor Guofang Wang for a lot of discussions on Sasakian geometry and useful suggestions during the preparation of this paper and thank professor Ildefonso  Castro for his comments and interests in this paper. Many thanks to professor Toru Sasahara for pointing out an error in Lemma \ref{main lem1}. The author is partially supported by the NSF of China(No.11501421).
{}
\vspace{1cm}\sc
Yong Luo

School of Mathematics and statistics,

Wuhan University, Wuhan 430072, China

{\tt yongluo@whu.edu.cn}

\vspace{1cm}\sc
\end{document}